\documentclass{article}\usepackage{amsmath,amssymb,latexsym,amsfonts,enumerate,theorem,longtable,amscd}
\newtheorem{definition}{Definition}\newtheorem{theorem}{Theorem}[section]\newtheorem{proposition}[theorem]{Proposition}
\newtheorem{corollary}[theorem]{Corollary}\newtheorem{lemma}[theorem]{Lemma}
\begin{document}

\begin{center}
{\large\bf Hopf algebra structure on symplectic superspace ${\rm SP}_q^{2|1}$}
\end{center}

\begin{center}
Salih Celik

Department of Mathematics, Yildiz Technical University, DAVUTPASA-Esenler, Istanbul, 34220 TURKEY.
\end{center}



\begin{abstract}
Super-Hopf algebra structure on the function algebra on the extended quantum symplectic superspace ${\rm SP}_q^{2|1}$ has been defined. The dual Hopf algebra is explicitly constructed.
\end{abstract}

\section{Introduction}
Quantum supergroups and quantum superalgebras are even richer mathematical subjects as compared to Lie supergroups and Lie superalgebras. A quantum superspace is a space that quantum supergroup acts with linear transformations and whose coordinates belong to a noncommutative associative superalgebra \cite{Manin2}.

Some algebras have been considered which are covariant with respect to the quantum supergroups in \cite{Chaic}. Using the corepresentation of the quantum supergroup ${\rm OSP}_q(1|2)$, some non-commutative spaces covariant under its coaction have been constructed \cite{Aiz-Chak1}. In the present work, we set up a super-Hopf algebra structure on an algebra which appears in both paper. We denoted this algebra by ${\mathbb O}({\rm SP}_q^{2|1})$. As is known, the matrix elements of the quantum supergroups ${\rm SP}_q(1|2)$ and ${\rm SP}_q(2|1)$ are the same and they act both quantum superspaces ${\rm SP}_q^{1|2}$ and
${\rm SP}_q^{2|1}$. But these two quantum superspaces are not the same. A study on ${\rm SP}_q^{1|2}$ was made in \cite{Celik}. Here we will work on the quantum symplectic superspace ${\rm SP}_q^{2|1}$.

\section{Review of quantum symplectic group}\label{sec2}

In this section, we will give some information about the structures of quantum symplectic groups as much as needed.

The algebra ${\mathbb O}({\rm OSP}_q(1|2))$ is generated by the {\it even} elements $a,b,c,d$ and {\it odd} elements $\alpha,\delta$. Standard FRT construction \cite{FRT} is obtained via the matrix $R$ given in \cite{Kul-Res}. Using the RTT-relations and the $q$-orthosymplectic condition, all defining relations of
${\mathbb O}({\rm OSP}_q(1|2))$ are explicitly obtained in \cite{Aiz-Chak1}:
\begin{theorem} 
The generators of ${\mathbb O}({\rm SP}_q(1|2))$  satisfy the relations
\begin{align} \label{2.1}
a b &= q^2 b a, \quad a c = q^2 c a, \quad a \alpha = q \alpha a, \quad a \delta = q \delta a + (q-q^{-1}) \alpha c, \nonumber\\
ad &= da + (q-q^{-1})[(1+q^{-1}) bc + q^{-1/2} \alpha\delta], \quad b c = c b, \quad bd = q^2 db, \nonumber\\
b \alpha &= q^{-1} \alpha b, \quad b \delta = q \delta b, \quad cd = q^2 dc, \quad c \alpha = q^{-1} \alpha c, \quad c \delta = q \delta c, \\
d \alpha &= q^{-1} \alpha d + (q^{-1}-q) \delta b, \quad d \delta = q^{-1} \delta d, \,\,\, \alpha\delta = -q \delta\alpha + q^{-1/2}(q^2-1) bc, \nonumber\\
\alpha^2 &= q^{1/2}(q-1) ba, \quad \delta^2 = q^{1/2}(q-1) dc.\nonumber
\end{align}
\end{theorem}

In (\ref{2.1}), the relations involving the elements $\gamma$, $e$ and $\beta$ are not written. They can be found in \cite{Aiz-Chak1}. Other relations that we need in this study are given below:
\begin{align} \label{2.2}
[e,\alpha]_q &= q^{1/2}(q-1) (\gamma b + \beta a), \quad [e,\beta]_{q^{-1}} = q^{-1/2}(q^{-1}-1) (\delta b + \alpha d), \nonumber\\
[e,\gamma]_q &= q^{1/2}(1-q) (\delta a + \alpha c), \quad [e,\delta]_{q^{-1}} = q^{-1/2}(1-q^{-1}) (\gamma d + \beta c), \nonumber\\
\beta^2 &= q^{1/2}(q-1) db, \quad \gamma^2 = q^{1/2}(q-1) ca, \\
e^2 &= 1-q^{-1/2}[\alpha,\delta]_q = 1+q^{1/2}[\beta,\gamma]_{q^{-1}} \nonumber
\end{align}
where $[u,v]_Q=uv-Qvu$.

The quantum superdeterminant is defined by
\begin{equation*}
D_q = ad - qbc - q^{1/2}\alpha\delta = da - q^{-1}bc + q^{-1/2}\delta\alpha.
\end{equation*}
This element of ${\mathbb O}({\rm SP}_q(2|1))$ commutes with all elements of ${\mathbb O}({\rm SP}_q(2|1))$.

If ${\mathbb A}$ and ${\mathbb B}$ are ${\mathbb Z}_2$-graded algebras, then their tensor product ${\mathbb A}\otimes{\mathbb B}$ is the ${\mathbb Z}_2$-graded algebra whose underlying space is ${\mathbb Z}_2$-graded tensor product of ${\mathbb A}$ and ${\mathbb B}$. The following definition gives the product rule for tensor product of algebras. Let us denote by $\tau(a)$ the {\it grade} (or {\it degree}) of an element $a\in {\mathbb A}$.

\begin{definition} 
If ${\mathbb A}$ is a ${\mathbb Z}_2$-graded algebra, then the product rule in the ${\mathbb Z}_2$-graded algebra ${\mathbb A}\otimes {\mathbb A}$ is defined by
\begin{equation*}
(a_1\otimes a_2)(a_3\otimes a_4) = (-1)^{\tau(a_2)\tau(a_3)} a_1a_3\otimes a_2a_4
\end{equation*}
where $a_i$'s are homogeneous elements in the algebra ${\mathbb A}$.
\end{definition}

\begin{definition}
A super-Hopf algebra is a vector space ${\mathbb A}$ over ${\mathbb K}$ with three linear maps
$\Delta:{\mathbb A}\longrightarrow {\mathbb A}\otimes{\mathbb A}$, called the coproduct, $\epsilon: {\mathbb A}\longrightarrow {\mathbb K}$, called the counit,
and $S:{\mathbb A}\longrightarrow{\mathbb A}$, called the coinverse, such that
\begin{align} \label{2.3}
(\Delta \otimes {\rm id}) \circ \Delta &= ({\rm id} \otimes \Delta) \circ \Delta, \\
m\circ(\epsilon\otimes{\rm id})\circ\Delta &= {\rm id} = m\circ({\rm id}\otimes\epsilon)\circ\Delta, \\
m\circ(S\otimes{\rm id})\circ\Delta &= \eta\circ\epsilon = m\circ({\rm id}\otimes S)\circ\Delta,
\end{align}
together with $\Delta({\bf 1}) ={\bf 1}\otimes{\bf 1}$, $\epsilon({\bf 1})=1$, $S({\bf 1})={\bf 1}$ and for any $a,b \in{\mathbb A}$
\begin{align} \label{2.6}
 \Delta(ab)=\Delta(a) \Delta(b), \quad \epsilon(ab)=\epsilon(a) \epsilon(b), \quad S(ab)=(-1)^{\tau(a)\tau(b)} S(b) S(a),
\end{align}
where $m:{\mathbb A} \otimes {\mathbb A} \longrightarrow {\mathbb A}$ is the product map, ${\rm id}:{\mathbb  A}\longrightarrow{\mathbb A}$ is the identity map and $\eta:{\mathbb K}\longrightarrow{\mathbb A}$.
\end{definition}

The super-Hopf algebra structure of ${\mathbb O}({\rm SP}_q(2|1))$ is given as usual in the following theorem.
\begin{theorem} 
There exists a unique super-Hopf algebra structure on the superalgebra ${\mathbb O}({\rm SP}_q(2|1))$ with co-maps $\Delta$, $\epsilon$ and $S$ such that
\begin{eqnarray*}
\Delta(t_{ij}) = \sum_{k=1}^3 t_{ik} \otimes t_{kj}, \quad \epsilon(t_{ij}) = \delta_{ij}, \quad S(T)=T^{-1}
\end{eqnarray*}
where
\begin{equation*}
T = \begin{pmatrix} a & \alpha & b \\ \gamma & e & \beta \\ c & \delta & d \end{pmatrix} = (t_{ij}), \quad
T^{-1} = \begin{pmatrix} d & q^{-1/2}\beta & -q^{-1}b \\ -q^{1/2}\delta & e & q^{-1/2}\alpha \\ -qc & -q^{1/2}\gamma & a \end{pmatrix}.
\end{equation*}
\end{theorem}

\section{Quantum symplectic superspace ${\rm SP}_q^{2|1}$}\label{sec3}
In this section, we define a super-Hopf algebra structure on the extended function algebra of the quantum superspace ${\rm SP}_q^{2|1}$.

\subsection{The algebra of polynomials on the quantum superspace ${\rm SP}_q^{2|1}$}\label{sec3.1}
The elements of the symplectic superspace are supervectors generated by two even and an odd components. We define a ${\mathbb Z}_2$-graded symplectic space
${\rm SP}^{2|1}$ by dividing the superspace ${\rm SP}^{2|1}$ of 3x1 matrices into two parts ${\rm SP}^{2|1}=V_0\oplus V_1$. A vector is an element of $V_0$ (resp. $V_1$) and is of grade 0 (resp. 1) if it has the form
\begin{equation*}
\begin{pmatrix} x \\ 0 \\ y\end{pmatrix}, \quad \mbox{resp.} \quad \begin{pmatrix} 0 \\ \theta \\ 0 \end{pmatrix}.
\end{equation*}
While the even elements commute to everyone, the odd element satisfies the relation $\theta^2=0$.

In \cite{Celik}, the quantum superspace ${\rm SP}_q^{2|1}$ is considered as the dual space of quantum superspace ${\rm SP}_q^{1|2}$  and then relations (\ref{3.1}) below are obtained by interpreting the coordinates as differentiations.

\begin{definition} 
Let ${\mathbb K}\langle x,\theta,y\rangle$ be a free associative algebra generated by $x$, $\theta$, $y$ and $I_q$ be a two-sided ideal generated by $x\theta-q\theta x$, $xy-q^2 yx$, $y\theta-q^{-1}\theta y$, $\theta^2-q^{1/2} (q-1) yx$. The quantum superspace ${\rm SP}_q^{2|1}$ with the function algebra
$${\mathbb O}({\rm SP}_q^{2|1}) = {\mathbb K}\langle x,\theta,y\rangle/I_q$$
is called ${\mathbb Z}_2$-graded quantum symplectic space (or quantum symplectic superspace).
\end{definition}
Here the coordinates $x$ and $y$ with respect to the ${\mathbb Z}_2$-grading are of grade 0 (or even), the coordinate $\theta$ with respect to the
${\mathbb Z}_2$-grading is of grade 1 (or odd).

According to the above definition, if $(x,\theta,y)^t\in {\rm SP}_q^{2|1}$ then we have
\begin{equation} \label{3.1}
x\theta = q\theta x, \quad \theta y = q y\theta, \quad yx=q^{-2} xy, \quad \theta^2=q^{1/2} (q-1) yx \end{equation}
where $q$ is a non-zero complex number. This associative algebra over the complex numbers is known as the algebra of polynomials over quantum (2+1)-superspace.

It is easy to see the existence of representations that satisfy (\ref{3.1}); for instance, there exists a representation
$\rho:{\mathbb O}({\rm SP}_q^{2|1})\to {\rm M}(3,{\mathbb C})$ such that matrices
\begin{equation*}
\rho(x) = \begin{pmatrix} q & 0 & 0 \\ 0 & q^2 & 0 \\ 0 & 0 & 1 \end{pmatrix}, \quad
\rho(\theta)=\begin{pmatrix} 0 & q-1 & 0 \\ 0 & 0 & 0 \\ q^{1/2} & 0 & 0 \end{pmatrix}, \quad
\rho(y) = \begin{pmatrix} 0 & 0 & 0 \\ 0 & 0 & 0 \\ 0 & 1 & 0 \end{pmatrix}
\end{equation*}
representing the coordinate functions satisfy the relations (\ref{3.1}).

Note that the last two relations in (\ref{3.1}) can be also written as a single relation. Therefore, we say that ${\mathbb O}({\rm SP}_q^{2|1})$ is the superalgebra with generators $x_\pm$ and $\theta$ satisfying the relations \cite{Chaic}
\begin{equation} \label{3.2}
x_\pm\theta = q^{\pm1}\theta x_\pm, \quad [x_+,x_-] = q^{-1/2} (q+1) \theta^2.
\end{equation}
where $x_+=x$ and $x_-=y$.

\begin{definition} \cite{Chaic} 
The quantum supersphere on the quantum symplectic superspace is defined by
\begin{equation*}
r = q^{1/2} x_- x_+ + \theta^2 - q^{-1/2} x_+ x_-.
\end{equation*}
\end{definition}

\subsection{A $\star$-structure on the algebra ${\mathbb O}({\rm SP}_q^{2|1})$}\label{sec3.2}

Here we define a ${\mathbb Z}_2$-graded involution on the algebra ${\mathbb O}({\rm SP}_q^{2|1})$.
\begin{definition} 
Let ${\mathbb A}$ be an associative superalgebra. A ${\mathbb Z}_2$-graded linear map $\star: {\mathbb A}\longrightarrow {\mathbb A}$ is called a superinvolution (or ${\mathbb Z}_2$-graded involution) if
\begin{equation*}
(ab)^\star=(-1)^{\tau(a)\tau(b)} b^\star a^\star, \qquad (a^\star)^\star=a
\end{equation*}
for any elements $a,b\in {\mathbb A}$. The pair $({\mathbb A},\star)$ is called a ${\mathbb Z}_2$-graded $\star$-algebra.
\end{definition}

If the parameter $q$ is real, then the algebra ${\mathbb O}({\rm SP}_q^{2|1})$ becomes a $\star$-algebra with involution determined by the following proposition.

\begin{proposition} 
If $q>0$ then the algebra ${\mathbb O}({\rm SP}_q^{2|1})$ supplied with the ${\mathbb Z}_2$-graded involution determined by
\begin{equation*} \label{3.3}
x^\star_+ = q^{1/2}\, x_-, \quad \theta^\star = {\bf i}\,\theta, \quad x^\star_- = q^{-1/2}\, x_+
\end{equation*}
becomes a super $\star$-algebra where ${\bf i}=\sqrt{-1}$.
\end{proposition}

\noindent{\it Proof}  We must show that the relations (\ref{3.2}) are invariant under the star operation. If $q$ is a positive number, we have
\begin{align*}
(x_\pm\theta-q^{\pm1}\theta x_\pm)^\star
& = ({\bf i}\,\theta)(q^{\pm1/2}x_\mp) - q^{\pm1} (q^{\pm1/2}x_\mp)({\bf i}\,\theta) \\
& = q^{\pm1/2}{\bf i}\, (\theta x_\mp - q^{\pm1} \theta x_\mp)
\end{align*}
and since $[x_+,x_-]^\star = [x_+,x_-]$
\begin{align*}
[x_+,x_-] & = [x_+,x_-]^\star = q^{-1/2}(q+1) (-\theta^\star \theta^\star) = q^{-1/2}(q+1) \theta^2.
\end{align*}
Hence the ideal $(x_\pm\theta-q^{\pm1}\theta x_\pm, \,[x_+,x_-]- q^{-1/2}(q+1) \theta^2)$ is $\star$-invariant and the quotient algebra
${\mathbb K}\langle x_+,\theta,x_-\rangle/(x_\pm\theta-q^{\pm1}\theta x_\pm, \,[x_+,x_-]- q^{-1/2}(q+1) \theta^2)$ becomes a $\star$-algebra. \hfill$\square$

\subsection{The super-Hopf algebra structure on ${\rm SP}_q^{2|1}$}\label{sec3.3}
We define the extended ${\mathbb Z}_2$-graded quantum symplectic space to be the algebra containing ${\rm SP}_q^{2|1}$, the unit and $x_+^{-1}$, the inverse of $x_+$, which obeys $x_+ x_+^{-1} = {\bf 1} =  x_+^{-1} x_+$. We will denote the unital extension of ${\mathbb O}({\rm SP}_q^{2|1})$ by
${\mathbb F}({\rm SP}_q^{2|1})$. The following theorem asserts that the superalgebra ${\mathbb F}({\rm SP}_q^{2|1})$ is a super-Hopf algebra:

\begin{theorem} 
The algebra ${\mathbb F}({\rm SP}_q^{2|1})$ is a ${\mathbb Z}_2$-graded Hopf algebra. The definitions of a coproduct, a counit and a coinverse on the algebra
${\mathbb F}({\rm SP}_q^{2|1})$ are as follows

\noindent(i) the coproduct $\Delta: {\mathbb F}({\rm SP}_q^{2|1})\longrightarrow {\mathbb F}({\rm SP}_q^{2|1})\otimes{\mathbb F}({\rm SP}_q^{2|1})$ is defined by
\begin{equation} \label{3.5}
\Delta(x_+) = x_+ \otimes x_+, \quad \Delta(\theta) = \theta\otimes{\bf 1}+{\bf 1}\otimes\theta, \quad \Delta(x_-) = x_+^{-1}\otimes x_- + x_-\otimes x_+^{-1}.
\end{equation}

\noindent(ii) the counit $\epsilon: {\mathbb F}({\rm SP}_q^{2|1})\longrightarrow {\mathbb C}$ is given by
\begin{equation*}
\epsilon(x_+)=1, \quad \epsilon(\theta)=0, \quad \epsilon(x_-)=0.
\end{equation*}

\noindent(iii) the algebra ${\mathbb F}({\rm SP}_q^{2|1})$ admits a ${\mathbb C}$-algebra antihomomorphism (coinverse)
$S: {\mathbb F}({\rm SP}_q^{2|1})\longrightarrow{\mathbb F}({\rm SP}_{q^{-1}}^{2|1})$ defined by
\begin{equation*}
S(x_+) = x_+^{-1}, \quad S(\theta) = - \theta, \quad S(x_-) = -x_+x_-x_+.
\end{equation*}
\end{theorem}

\noindent{\it Proof} The axioms (\ref{2.3})-(5) are satisfied automatically. It is also not difficult to show that the co-maps preserve the relations (\ref{3.2}). In fact, for instance,
\begin{align*}
\Delta([x_+,x_-])
& = \Delta(x_+ x_- - x_- x_+) = {\bf 1}\otimes[x_+,x_-] + [x_+,x_-]\otimes{\bf 1}\\
& = q^{-1/2} (q+1) ({\bf 1}\otimes\theta^2 + \theta^2\otimes{\bf 1})\\
\Delta(\theta^2)
& = {\bf 1}\otimes\theta^2 + \theta^2\otimes{\bf 1},
\end{align*}
and
\begin{align*}
S([x_+,x_-])
& = -[x_+, x_-], \quad S(\theta^2) = -\theta^2.
\end{align*}
Since $S^2(a)={\rm id}(a)$ for all $a\in {\mathbb F}({\rm SP}_q^{2|1})$, the coinverse $S$ is of second order. \hfill$\square$

The set $\{x^k\theta^ly^m: \, k,l,m\in{\mathbb N}_0\}$ form a vector space basis of ${\mathbb F}({\rm SP}_q^{2|1})$. The formula (\ref{3.5}) gives the action of the coproduct $\Delta$ only on the generators. The action of $\Delta$ on product on generators can be calculated by taking into account that $\Delta$ is a homomorphism.

\begin{corollary} 
For the quantum supersphere $r$, we have
\begin{equation*}
\Delta(r) = r \otimes {\bf 1} + {\bf 1} \otimes r, \quad \epsilon(r) = 0, \quad S(r) = - r.
\end{equation*}
\end{corollary}

\noindent{\it Proof} Using the definition of $\Delta$, as an algebra homomorphism, on the generators of ${\mathbb F}({\rm SP}_q^{2|1})$ in (\ref{3.5}), it is easy to see that the element $r\in {\mathbb F}({\rm SP}_q^{2|1})$ is a primitive element, that is,
\begin{align*}
\Delta(r)
& = q^{1/2}(x_+^{-1}\otimes x_- + x_- \otimes x_+^{-1})(x_+ \otimes x_+) + (\theta\otimes{\bf 1} + {\bf 1}\otimes\theta)(\theta\otimes{\bf 1} + {\bf 1}\otimes\theta)\\
& \quad - q^{-1/2} (x_+ \otimes x_+)(x_+^{-1}\otimes x_- + x_- \otimes x_+^{-1})\\
& = q^{1/2} ({\bf 1}\otimes x_- x_+ + x_-x_+ \otimes {\bf 1}) + \theta^2\otimes{\bf 1} + {\bf 1}\otimes\theta^2 \\ & \quad - q^{-1/2} ({\bf 1}\otimes x_+ x_- + x_+x_- \otimes {\bf 1})\\
& = {\bf 1}\otimes (q^{1/2} x_- x_+ + \theta^2 - q^{-1/2} x_+ x_-) + (q^{1/2} x_- x_+ + \theta^2 - q^{-1/2} x_+ x_-)\otimes{\bf 1}.
\end{align*}
Since $\epsilon({\bf 1})=1$ and
\begin{equation*}
m({\rm id}\otimes\epsilon)\Delta(r) = r \epsilon({\bf 1}) + \epsilon(r) {\bf 1} = r = m(\epsilon\otimes{\rm id})\Delta(r),
\end{equation*}
we obtain $\epsilon(r)=0$. Finally, using the fact that $S$ is an anti-homomorphism we get
\begin{align*}
S(r)
& = q^{1/2} x_+^{-1} (-x_+ x_- x_+) - (-\theta)(-\theta) - q^{-1/2} (-x_+ x_- x_+)x_+^{-1}\\
& = -(q^{1/2} x_- x_+ + \theta^2 - q^{-1/2} x_+ x_-),
\end{align*}
as desired. \hfill$\square$

\subsection{Coactions on the quantum symplectic superspace}\label{sec3.4}
Let $a,b,c,d,e,\gamma,\alpha,\delta,\beta$ be elements of an algebra ${\mathbb A}$. Assuming that the generators of ${\mathbb O}({\rm SP}_q(2|1))$ super-commute with the elements of ${\mathbb O}({\rm SP}_q^{2|1})$, define the components of the vectors $X'=(x', \theta', y')^t$ and $X''=(x'', \theta'', y'')^t$ using the following matrix equalities
\begin{equation} \label{3.7}
X' = T\, X \quad {\rm and} \quad (X'')^t = X^t \,T
\end{equation}
where $X=(x, \theta, y)^t \in {\rm SP}_q^{2|1}$. If we assume that $q\ne1$ then we have the following theorem proving straightforward computations.

\begin{theorem} 
If the couples $(x', \theta', y')$ and $(x'', \theta'', y'')$ in $(\ref{3.7})$ satisfy the relations $(\ref{3.1})$, then the generators of
${\mathbb O}({\rm SP}_q(2|1))$ together with $q$-orthosymplectic condition fulfill the relations $(\ref{2.1})$.
\end{theorem}

A left quantum space (or left comodule algebra) for a Hopf algebra $H$ is an algebra ${\mathbb X}$ together with an algebra homomorphism (left coaction)
$\delta_L:{\mathbb X}\longrightarrow H\otimes{\mathbb X}$ such that
\begin{equation*}
(\mbox{id}\otimes\delta_L)\circ\delta_L=(\Delta\otimes\mbox{id})\circ\delta_L \quad \mbox{and} \quad (\epsilon\otimes\mbox{id})\circ\delta_L=\mbox{id}.
\end{equation*}
Right comodule algebra can be defined in a similar way.

\begin{theorem} 
(i) The algebra ${\mathbb O}({\rm SP}_q^{2|1})$ is a left and right comodule algebra of the Hopf algebra ${\mathbb O}({\rm SP}_q(2|1))$ with left coaction
$\delta_L$ and right coaction $\delta_R$ such that
\begin{equation} \label{3.8}
\delta_L(X_i) = \sum_{k=1}^3 t_{ik} \otimes X_k, \quad \delta_R(X_i) = \sum_{k=1}^3 X_k \otimes t_{ki}.
\end{equation}

\noindent(ii) The quantum supersphere $r$ belongs to the center of ${\mathbb O}({\rm SP}_q^{2|1})$ and satisfies $\delta_L(r)={\bf 1}\otimes r$ and $\delta_R(r)=r \otimes{\bf 1}$.
\end{theorem}

\noindent{\it Proof}
(i): These assertions are obtained from the relations in (\ref{2.1}) and (\ref{2.2}) together with (\ref{3.1}).

(ii): That $r$ is a central element of ${\mathbb O}({\rm SP}_q^{2|1})$ is shown by using the relations in (\ref{3.1}). To show that $\delta_L(r)={\bf 1}\otimes r$ and $\delta_R(r)=r \otimes{\bf 1}$ we use the definitions of $\delta_L$ and $\delta_R$ in (\ref{3.8}) and the relations (\ref{2.1}) and (\ref{2.2}) with $D_q={\bf 1}$. \hfill$\square$

\section{Dual pairings of Hopf algebras}\label{sec4}
Let us consider the algebra ${\mathbb A}$ as a Hopf algebra. Then its dual ${\mathbb U}\doteq {\mathbb A}'$ is a Hopf algebra as well. Using the coproduct $\Delta$ in ${\mathbb A}$, one defines a product in ${\mathbb U}$. Using the product in the Hopf algebra ${\mathbb A}$, one defines a coproduct in ${\mathbb U}$.

In this section, in order to obtain the dual of the Hopf algebra ${\mathbb A}$ defined in Section 3, we have applied to
${\mathbb A}={\mathbb F}({\rm SP}_q^{2|1})$ the approach given in \cite{Sudb}. We are going to use the relations (\ref{3.1}) in order to facilitate the
operations in this section. Let us begin with the definition of the duality \cite{Abe}.
\begin{definition} 
A dual pairing of two super bialgebras ${\mathbb U}$ and ${\mathbb A}$ is a bilinear mapping
\begin{align} \label{4.1}
\langle,\rangle:{\mathbb U}\times{\mathbb A}\longrightarrow {\mathbb K}, \quad (u,a)\mapsto \langle u,a\rangle, \quad  u \in {\mathbb U}, \, a \in {\mathbb A}
\end{align}
such that
\begin{align} \label{4.2}
\langle u, ab\rangle
&= \langle\Delta_{\mathbb U}(u), a \otimes b\rangle, \quad \langle uv, a\rangle = \langle u \otimes v, \Delta_{\mathbb A}(a)\rangle, \\
\langle u, 1_{\mathbb A}\rangle
&= \epsilon_{\mathbb U}(u), \quad \langle 1_{\mathbb U}, a\rangle = \epsilon_{\mathbb A}(a) \nonumber
\end{align}
for all $u, v \in {\mathbb U}$ and $a, b \in {\mathbb A}$.
\end{definition}
We say that the pairing is non-degenerate if
$$\langle u,a\rangle = 0 \quad (\forall a \in {\mathbb A}) \Longrightarrow u = 0 \quad {\rm and} \quad
  \langle u,a\rangle = 0 \quad (\forall u \in {\mathbb U}) \Longrightarrow a = 0.$$

Two super Hopf algebras ${\mathbb U}$ and ${\mathbb A}$ are said to be in duality if they are in duality as bialgebras and if
\begin{align*}
\langle S_{\mathbb U}(u), a\rangle = \langle u, S_{\mathbb A}(a)\rangle, \quad \forall u \in {\mathbb U}, \quad a \in {\mathbb A}.
\end{align*}

Such a pairing can be extended to a pairing of ${\mathbb U} \otimes {\mathbb U}$ and ${\mathbb A} \otimes {\mathbb A}$ by
\begin{align*}
\langle u \otimes v, a \otimes b\rangle = (-1)^{\tau(v)\tau(a)}\langle u,a\rangle \langle v,b\rangle.
\end{align*}

It is enough to define the pairing (\ref{4.1}) between the generating elements of the two algebras. Pairing for any other elements of ${\mathbb U}$ and ${\mathbb A}$ follows from (\ref{4.2}) and the bilinear form inherited by the tensor product.

As a super Hopf algebra ${\mathbb A}={\mathbb F}({\rm SP}^{2|1}_q)$ is generated by the elements $x$, $\theta$, $y$ and a basis is given by all monomials of the form
$$f = x^k \theta^l y^n$$
where $k,l,n \in {\mathbb N}_0$. Let us denote the dual algebra by ${\mathbb U}$ and its generating elements by $K$, $\nabla$ and $N$.

\begin{theorem} 
The commutation relations between the generators of the algebra ${\mathbb U}$ dual to ${\mathbb A}$ are the following:
\begin{align*}
K\nabla = \nabla K, \quad KN = NK, \quad \nabla N = N\nabla, \quad \nabla^2=0.
\end{align*}
\end{theorem}

\noindent{\it Proof}
The pairing is defined through the tangent vectors as follows
\begin{equation*}
\langle K, f\rangle = k \delta_{l,0}\delta_{n,0}, \quad \langle\nabla, f\rangle = \delta_{l,1}, \quad \langle N, f\rangle = \delta_{n,1}.
\end{equation*}
We also have
\begin{align*}
\langle 1_{\mathbb U}, f\rangle = \epsilon_{\mathbb A}(f) = \delta_{k,0}.
\end{align*}
Using the defining relations one gets
\begin{align*}
\langle KN, f\rangle = (k-1)(\delta_{l,0}\delta_{n,1}+\delta_{l,1}\delta_{n,1}) = \langle NK,f\rangle,
\end{align*}
where differentiation is from the right as this is most suitable for differentiation in this basis. Thus one obtains
\begin{align*}
\langle KN-NK, f\rangle = 0.
\end{align*}
The other relations can be obtained similarly. \hfill$\square$

\begin{theorem} 
The Hopf algebra structure of the algebra ${\mathbb U}$ is given by:

\noindent(i) the comultiplication is given by
\begin{align*}
\Delta_{\mathbb U}(K) &= K\otimes 1_{\mathbb U} + 1_{\mathbb U} \otimes K, \nonumber\\
\Delta_{\mathbb U}(\nabla) &= \nabla \otimes 1_{\mathbb U} + q^K \otimes \nabla, \nonumber\\
\Delta_{\mathbb U}(N) &= N \otimes q^{-2K} + 1_{\mathbb U} \otimes N,
\end{align*}
(ii) the counity is given by
\begin{align*}
\epsilon_{\mathbb U}(K) = 0, \quad \epsilon_{\mathbb U}(\nabla) = 0, \quad \epsilon_{\mathbb U}(N) = 0,
\end{align*}
(iii) the coinverse is given by
\begin{align*}
S_{\mathbb U}(K) = - K, \quad S_{\mathbb U}(\nabla) = - q^{-K}\nabla, \quad S_{\mathbb U}(N) = - q^{2K}N.
\end{align*}
\end{theorem}

\noindent{\it Proof}
We only will obtain the actions of $\Delta_{\mathbb U}$, $\epsilon_{\mathbb U}$ and $\kappa_{\mathbb U}$ on $\nabla$. The others can be obtained similarly. So, we assume that the action of $\Delta_{\mathbb U}$ on $\nabla$ is $\Delta_{\mathbb U}(\nabla)= \nabla \otimes c_1 + c_2 \otimes \nabla$. Then the commutation relations in ${\mathbb A}$ (or comultiplication in ${\mathbb U}$) will imply that $c_1=q^{-K}c_2$. Indeed,
\begin{align*}
\langle\Delta_{\mathbb U}(\nabla), \theta\otimes x^k\rangle
&= \langle\nabla \otimes c_1 + c_2 \otimes \nabla, \theta\otimes x^k\rangle \\
&= \langle\nabla,\theta\rangle \langle c_1,x^k\rangle - \langle c_2,\theta\rangle \langle \nabla,x^k\rangle = \langle c_1,x^k\rangle
\end{align*}
and
\begin{align*}
\langle\Delta_{\mathbb U}(\nabla), x^k\otimes \theta\rangle
&= \langle\nabla \otimes c_1 + c_2 \otimes \nabla, x^k\otimes \theta\rangle\\
&= \langle\nabla,x^k\rangle \langle c_1,\theta\rangle + \langle c_2,x^k\rangle \langle\nabla,\theta\rangle = \langle c_2,x^k\rangle
\end{align*}
so that, considering $\langle\nabla,x^k\theta\rangle$ and $\langle \nabla,\theta x^k\rangle$ and taking $c_1=1_{\mathbb U}$ we have found the action of
$\Delta_{\mathbb U}$ on $\nabla$.

The action of $\epsilon_{\mathbb U}$ on $\nabla$ is
\begin{align*}
\epsilon_{\mathbb U}(\nabla) = \langle \nabla,1_{\mathbb A}\rangle = \langle \nabla,x^0\rangle = 0.
\end{align*}
This can also be found from the identities
\begin{equation*}
\mu \circ (\epsilon_{\mathbb U} \otimes 1_{\mathbb U})\Delta_{\mathbb U}(\nabla) = 1_{\mathbb U}(\nabla) = \mu \circ (1_{\mathbb U} \otimes \epsilon_{\mathbb U})\Delta_{\mathbb U}(\nabla).
\end{equation*}

The action of $S_{\mathbb U}$ on $\nabla$ is
\begin{align*}
\langle S_{\mathbb U}(\nabla),\theta x^k \rangle
&= \langle \nabla, S_{\mathbb A}(\theta x^k) \rangle = \langle \nabla, x^{-k}(-\theta) \rangle = \langle \nabla, -q^{-k} \theta x^{-k} \rangle = -q^{-k} \\
&\Longrightarrow  S_{\mathbb U}(\nabla) = - q^{-K} \nabla.
\end{align*}
This can also be found from the identities
\begin{equation*}
\mu \circ (S_{\mathbb U} \otimes 1_{\mathbb U})\Delta_{\mathbb U}(\nabla) = \epsilon_{\mathbb U}(\nabla) = \mu \circ (1_{\mathbb U} \otimes S_{\mathbb U})\Delta_{\mathbb U}(\nabla). \hspace*{4cm} \square
\end{equation*}

\section{An $h$-deformation of the superspace ${\rm SP}^{2|1}$} \label{sec5}

In this section, we introduce an $h$-deformation of the superspace ${\rm SP}^{2|1}$ from the $q$-deformation via a contraction following the method of \cite{agha}. Consider the $q$-deformed algebra of functions on the quantum superspace ${\rm SP}_q^{2|1}$ generated by $x_\pm$ and $\theta$ with the relations (\ref{3.2}).

We introduce new coordinates $X_\pm$ and $\Theta$ by
\begin{equation*} \label{5.1}
{\bf x} = \begin{pmatrix} x_+ \\ \theta \\ x_- \end{pmatrix} = \begin{pmatrix} 1 & 0 & 0 \\ 0 & 1 & 0 \\ \frac{h}{q-1} & 0 & 1 \end{pmatrix} \begin{pmatrix} X_+ \\ \Theta \\ X_- \end{pmatrix} = g \,{\bf X}.
\end{equation*}
When the relations (\ref{3.2}) are used, taking the limit $q\to1$ we obtain the following exchange relations, which define the $h$-superspace
${\rm SP}_h^{2|1}$:
\begin{definition} 
Let ${\mathbb O}({\rm SP}_h^{2|1})$ be the algebra with the generators $X_\pm$ and $\Theta$ satisfying the relations
\begin{equation} \label{5.2}
X_+\Theta = \Theta X_+, \quad X_- \Theta = \Theta X_- - 2h \Theta X_+, \quad X_+ X_- = X_- X_+ + 2\Theta^2
\end{equation}
where the coordinates $X_\pm$ are even and the coordinate $\Theta$ is odd. We call ${\mathbb O}({\rm SP}_h^{2|1})$ the algebra of functions on the
${\mathbb Z}_2$-graded quantum space ${\rm SP}_h^{2|1}$.
\end{definition}

$h$-deformed supersphere on the symplectic $h$-superspace is given by
\begin{equation*} \label{5.3}
r_h = X_- X_+ + \Theta^2 + h \,X_+^2 - X_+ X_- = h \,X_+^2 - \Theta^2.
\end{equation*}

It is easily seen that the quantum supersphere $r_h$ belongs to the center of the superalgebra ${\mathbb O}({\rm SP}_h^{2|1})$.

The definition of dual $q$-deformed symplectic superspace is given as follows \cite{Aiz-Chak1}.
\begin{definition} 
Let ${\mathbb K}\{\varphi_+,z,\varphi_-\}$ be a free associative algebra generated by $z$, $\varphi_+$, $\varphi_-$ and $I_q$ be a two-sided ideal generated by $z\varphi_\pm-q^{\pm1} \varphi_\pm z$, $\varphi_- \varphi_+ + q^{-2}\varphi_+ \varphi_- + q^{-2}Qz^2$ and $\varphi_\pm^2$. The quantum superspace $SP_q^{1|2}$ with the function algebra
$${\mathbb O}(SP_q^{1|2}) = {\mathbb K}\{\varphi_+,z,\varphi_-\}/I_q$$
is called ${\mathbb Z}_2$-graded quantum symplectic space (or quantum symplectic superspace) where  $Q=q^{1/2}-q^{3/2}$ and $q\ne0$.
\end{definition}
In case of exterior $h$-superspace, we use the transformation
\begin{equation*}
\hat{\bf x}=g\hat{\bf X}
\end{equation*}
with the components $\varphi_+$, $z$ and $\varphi_-$ of $\hat{\bf x}$. The definition is given below.

\begin{definition} 
Let $\Lambda({\rm SP}_h^{2|1})$ be the algebra with the generators $\Phi_\pm$ and $Z$ satisfying the relations
\begin{align*} \label{5.4}
\Phi_+ Z &= Z \Phi_+, \quad Z \Phi_- = \Phi_- Z - 2h \Phi_+ Z, \quad \Phi_- \Phi_+ = - \Phi_+ \Phi_-, \nonumber\\
\Phi_+^2 &= 0, \quad \Phi_-^2 = h(2\Phi_- \Phi_+ - Z^2)
\end{align*}
where the coordinate $Z$ is even and the coordinates $\Phi_\pm$ are odd. We call $\Lambda({\rm SP}_h^{2|1})$ the quantum exterior algebra of the
${\mathbb Z}_2$-graded quantum space ${\rm SP}_h^{2|1}$.
\end{definition}

\section{A new algebra derived from ${\mathbb F}({\rm SP}_q^{2|1})$}\label{sec6}


It is known that an element of a Lie group can be represented by exponential of an element of its Lie algebra. By virtue of this fact, one can define the
generators of the algebra ${\mathbb F}({\rm SP}_q^{2|1})$ as
\begin{equation} \label{6.1}
x_+ := e^u, \quad \theta := q^{-1/2} \,\xi, \quad x_- := e^{-u} v.
\end{equation}
Then, the following lemma can be proved by direct calculations using the relations
\begin{equation*}
x_\pm^k \theta = q^{\pm k} \,\theta x_\pm^k, \quad
[x_+^k,x_-] = q^{-1/2} \,\frac{q^{2k}-1}{q-1} \,\theta^2 \,x_+^{k-1}, \quad \forall k\ge1
\end{equation*}
whose the proof follows from induction on $k$.

\begin{lemma} 
The generators $u$, $\xi$, $v$ have the following commutation relations (Lie (anti-)brackets)
\begin{equation} \label{6.2}
[u,\xi] = \hbar \,\xi, \quad [\xi,v] = 0, \quad [u,v] = \frac{2\hbar}{1 - e^{-\hbar}} \,\xi^2,
\end{equation}
where $q=e^\hbar$ and $\hbar\in{\mathbb R}$.
\end{lemma}

We denote the algebra for which the generators obey the relations (\ref{6.2}) by ${\mathbb L}_{\hbar}:={\mathbb L}({\rm SP}_q^{2|1})$. The ${\mathbb Z}_2$-graded Hopf algebra structure of ${\mathbb L}_{\hbar}$ can be read off from Theorem 3.2:

\begin{theorem} 
The algebra ${\mathbb L}_{\hbar}$ is a ${\mathbb Z}_2$-graded Hopf algebra. The definitions of a coproduct, a counit and a coinverse on the algebra
${\mathbb L}_{\hbar}$ are as follows:
\begin{equation*} \label{6.3}
\Delta(u_i) = u_i\otimes{\bf 1} + {\bf 1}\otimes u_i, \quad \epsilon(u_i) = 0, \quad S(u_i) = - u_i
\end{equation*}
for $u_i\in\{u,\xi,v\}$.
\end{theorem}

The following proposition can be easily proved by using the Proposition 3.1 together with (\ref{6.1}).
\begin{proposition} 
The algebra ${\mathbb L}_{\hbar}$ supplied with the ${\mathbb Z}_2$-graded involution determined by
\begin{equation*} \label{6.4}
u^\star = \tfrac{1}{2}\, \hbar + \ln(e^{-u} v), \quad \xi^\star = {\bf i}\,\xi, \quad v^\star = v
\end{equation*}
becomes a super Lie $\star$-algebra.
\end{proposition}

\end{document}